\documentclass[11pt]{article}
\usepackage[utf8]{inputenc}
\usepackage{amsfonts,epsf,amsmath,amssymb,here}
\usepackage[numbers]{natbib}
\usepackage{epsfig}
\usepackage{latexsym,multirow}
\usepackage{pgf,tikz,subfigure}
\usetikzlibrary{decorations.markings}
\usepackage{algorithm}
\usepackage[noend]{algorithmic}
\usepackage{paralist}
\usepackage{url}

\newtheorem{theorem}{\bf Theorem}[section]

\newtheorem{lemma}[theorem]{\bf Lemma}

\newtheorem{definition}[theorem]{\bf Definition}

\newcommand{\proof}{\noindent{\bf Proof.\ }}
\newcommand{\qed}{\hfill $\square$ \bigskip}

\DeclareMathOperator*{\id}{id}
\DeclareMathOperator*{\init}{init}

\newcommand{\cT}{{\cal T}}
\begin{document}

\title{An algebraic approach to enumerating non-equivalent double traces in graphs}

\author{
Nino Ba\v{s}i\'{c} \\
Faculty of Mathematics and Physics, University of Ljubljana \\
Jadranska 19, 1000 Ljubljana, Slovenia  \\
\url{nino.basic@fmf.uni-lj.si} 
\and
Drago Bokal \\
Faculty of Natural Sciences and Mathematics, University of Maribor \\
Koro\v ska 160, 2000 Maribor, Slovenia \\
\url{bokal@uni-mb.si}
\and
Tomas Boothby\\
Department of Mathematics, Simon Fraser University \\
Burnaby, B.C. V5A 1S6, Canada \\
\url{tboothby@sfu.ca}
\and 
Jernej Rus \\
IMFM \\
Jadranska 19, 1000 Ljubljana, Slovenia \\
\url{jernej.rus@gmail.com}
}
\date{\today}
\maketitle

\begin{abstract}
Recently designed biomolecular approaches to build single chain 
polypeptide polyhedra as molecular origami nanostructures have 
risen high interest in various double traces of the underlying graphs of 
these polyhedra. Double traces are walks that traverse every edge of 
the graph twice, usually with some additional conditions on 
traversal direction and vertex neighborhood coverage. Given that 
double trace properties are intimately related to the efficiency of 
polypeptide polyhedron construction, enumerating all different possible 
double traces and analyzing their properties is an important step in the 
construction. In the paper, we study the automorphism group of double traces and present an algebraic approach to this problem, yielding a branch-and-bound 
algorithm.
\end{abstract}


\medskip\noindent
{\bf Math. Subj. Class. (2010):}
05C30, 
05C45, 
05C85, 
68R10, 
92E10. 

\section{Introduction}
\label{sec:intro}

Gradi\v sar et al.\ presented a novel self-assembly strategy for 
polypeptide nanostructure design in $2013$~\cite{gr-2013}. Their 
research was already improved by Kočar et al.\ in 2015, who 
developed another alternative strategy to design topofolds --- 
nanostructures built from polypeptide arrays of interacting modules 
that define their topology~\cite{ko-2015}. Such approaches are paving 
the way to a significant breakthrough in the field of protein origami, 
an area advancing a step ahead of DNA origami, where many researchers 
have spent the better part of the past decade by folding the molecules 
into dozens of intricate shapes. 

A polyhedron $P$ that is composed from a single polymer chain can 
be naturally represented by a graph $G(P)$ of the polyhedron. As 
every edge of $G(P)$ corresponds to a coiled-coil dimer in the 
self-assembly process, exactly two biomolecular segments are associated 
with every edge of $G(P)$. Hence, every edge of $G(P)$ is in its 
biomolecular structure replaced by two copies, resulting in a graph 
$G'(P)$ obtained from $G(P)$ by replacing every edge with a digon.
The graph $G'(P)$ is therefore Eulerian, and its Eulerian walks 
(i.e., walks that traverse every edge of $G(P)$ precisely twice),
called \textit{double traces} of $G(P)$, play a key role in modeling the 
construction process. Note that the argument shows that every graph
admits a double trace.

Double traces with additional properties related to stability of the 
constructed polyhedra were introduced as a combinatorial model underlying these 
approaches to polypeptide polyhedra design in~\cite{kl-2013} and~\cite{fi-2013}. 
Stability of the resulting polyhedron depends on two additional properties:
one relates to whether in the double trace the neighborhoods of vertices can be split,
and the other defines whether the edges of the double trace are traversed
twice in the same or in different directions. 

To define the first property, let an alternate sequence
$W=w_0 e_1 w_1\ldots w_{2m-1} e_{2m} w_{2m}$, where $e_i$ is an edge between vertices $w_{i-1}$ and $w_i$, be a double trace --- a closed walk which traverses every edge of graph exactly twice. Note that we always consider vertex sequence of a double trace with indices taken modulo $2m$. 
(Since the graph $G(P)$ is simple, so are all our
other graphs, except $G'(P)$. Hence, a double trace is completely 
described by listing the vertices of the corresponding walk and we sometimes write double trace as a sequence consisting only from vertices.) 
For a set of vertices $N \subseteq N(v)$, a double trace $W$ has a \textit{$N$-repetition} at vertex $v$ (nontrivial $N$-repetition in~\cite{fi-2013}),
if $N$ is nonempty, $N \neq N(v)$, and whenever $W$ comes to $v$ from a vertex in $N$ it also continues to a vertex in $N$. More formally $W$ has a $N$ repetition at $v$ if the following implication holds:
\begin{equation*}
\text{\emph{for every $i \in \{0,\ldots,\ell-1\}$: if $v=w_i$ then $w_{i+1} \in N$ if and only if $w_{i-1} \in N$.}}
\label{eq:repetition}
\end{equation*}
Then, $W$ is a \textit{strong trace} 
if $W$ is for every veretex $v$ without $N$-repetitions at $v$.  
It is a nontrivial result of~\cite{fi-2013} that every graph admits a strong trace.
A weaker concept of \textit{$d$-stable trace} requires that whenever $W$
has an $N$-repetition at some vertex $v$, then $|N|>d$. Fijavž et al.\ showed that
$G$ admits a $d$-stable trace if and only if $\delta(G)\ge d$~\cite{fi-2013}.

For the second property, note that there are precisely two directions 
to traverse an edge $e=uv$. If the same direction is used both times 
$W$ traverses $e$, then $e$ is a \textit{parallel} edge w.r.t.\ $W$, 
otherwise it is an \textit{antiparallel} edge. A double trace $W$ is 
\textit{parallel}, if all edges of $G$ are parallel w.r.t.\ $W$ and
is \textit{antiparallel}, if all the edges are antiparallel.
Interestingly, antiparallel traces appeared (under a different name)
two centuries ago in a study of properties of labyrinths by Tarry~\cite{tarry-1895},
who observed (in our language) that every connected graph admits an 
antiparallel double trace. Fijavž et al.\ extended this by characterizing
the graphs that admit an antiparallel strong trace~\cite{fi-2013},
and Rus upgraded the result to characterize graphs that admit 
an antiparallel $d$-stable trace~\cite{rus-2015}. The former characterization can be 
algorithmically verified using algorithms of~\cite{fur-1988}, but regarding
the latter, it is only known that the existence of antiparallel $1$-stable traces can be verified
using Thomassen's modification of the aforementioned algorithm, as published
in~\cite{th-1990} and later corrected by Benevant L\' opez and Soler Fern\' andez in~\cite{ben-1998}.
Similar modification of algorithm for spanning tree parity problem presented in~\cite{gab-1986} would work for $d > 1$ as well,
rendering the problem ``Does there exist an antiparallel $d$-stable trace in $G$?"
polynomially tractable. Some additional research was also made in~\cite{br-1988} and~\cite{el-2013}.

It is easy to obtain new traces from a given trace: one can change direction
of tracing or start at a different vertex. Also, if graph possesses certain symmetries, 
these may reflect in the trace. Such changes do not alter any properties of the
trace, hence we call the resulting traces equivalent, and we are interested in
non-equivalent traces, as introduced in~\cite{kl-2013}: 

\begin{definition}
\label{def:equal}
Two double traces $W$ and $W'$ are called {\em equivalent} if $W'$ can be obtained from $W$ (i) by reversion of $W,$ (ii) by shifting $W,$ 
(iii) by applying a permutation on $W$ induced by an automorphism of $G$,
or (iv) using any combination of the previous three operations.
If that is not the case, $W$ and $W'$ are {\em non-equivalent}.
\end{definition}
Two double traces $W$ and $W'$ are called {\em different} if their vertex 
sequences are not the same. Two different double traces may be equivalent.
It is easy to see that equivalence of double traces is an equivalence relation on
the {\em set $\cT$ of all different double traces}, and hence on any subset (such as 
strong traces, $d$-stable traces etc.). The main contribution of our paper
is designing for each of the subsets of interest an algorithm that, 
for a given graph as an input, outputs
precisely one representative of each equivalence class.
This representative will be the unique minimal element for 
the following linear ordering, called \textit{lexicographical ordering} 
of double traces. We assume that the vertices of $G$ are linearly ordered
as $v_0<v_1<\ldots<v_{n-1}$, and that $v_0$, $v_1$ are adjacent.
This linear ordering induces an ordering on the set of double traces
as follows:

\begin{definition}
\label{def:lex}
Given two double traces $W = w_0 \ldots w_{2m}$ and $W' = w_0' \ldots w_{2m}'$, $W$ is {\em lexicographically smaller or equal} to $W'$, denoted $W \le_{\mathit{lex}} W'$, if and only if 
$W=W'$ or the first $w_i$, which is different from $w_i'$, is smaller than $w_i'$.
\end{definition}

As lexicographical order is a linear order, it is clear that 
any finite set $S$ of double
traces has a unique lexicographically smallest member. We call that member the 
\textit{canonical representative} of $S$. 

For a more detailed treatment of double-trace related definitions we 
refer the reader to~\cite{fi-2013}. For other terms and concepts from 
graph theory not defined here, we refer to~\cite{we-1996}. 

Let the {\em automorphism group} $\text{Aut}(G)$ of $G$ be denoted by $A$. 
An automorphism $\pi\in A$ acts on $\cT$ by mapping a double trace 
$W = w_0 \ldots w_{2m}$ to $\pi(W) = \pi(w_0) \ldots \pi(w_{2m})$.
Let $\rho:\cT\rightarrow\cT$ be a reversal that maps $W = w_0 \ldots w_{2m}$
to $W' = w_{0} w_{2m} \ldots w_{1}$, and, for $i=0,\ldots,2m$, 
let $\sigma_i$ be an $i$-shift that maps
$W = w_0 \ldots w_{2m}$ to $W'' = w_i \ldots w_{2m+i}$.
Note that $\sigma_0=\sigma_{2m}=\id$.
Then the group $A$, the group $R=\{\id,\rho\}$, and the group $S=\{\sigma_i \mid i=0,\ldots,2m-1\}$ are
three groups acting on $\cT$ (or any of its subsets). 
Note that groups $R$ and $S$ do not commute and $\langle R,S \rangle$ is a dihedral group of symmetries of a regular $(2m)$-gon, where $E(G) = 2m$.
Therefore the orbits of the direct product $\Gamma=A \times \langle R, S \rangle$
are precisely the equivalence classes of double traces for the relation from Definition~\ref{def:equal}. Hence, a canonical representative of each
equivalence class is the lexicographically smallest element of each class.
We say that a double trace is \textit{canonical}, if it is the 
lexicographically smallest element of its orbit, meaning that every element
of $\Gamma$ maps it to a lexicographically larger (or equal) element.
Note that to verify canonicity of a particular double trace,
it is not enough to check whether the generators of $\Gamma$ 
map it to a larger element (we leave finding an example to the reader). 

It is easy to see that every canonical double trace
starts with $v_0v_1$ (by assumption, these two vertices are adjacent) and that every double trace is equivalent to at least one canonical double trace.
Double traces (not necessary canonical) starting with $v_0v_1$ are called {\em simple}.
More details on graph automorphisms can be found in~\cite{god-2001},
but we do conclude this introduction with an example of the action of 
$\Gamma$ on $\cT$ in the case of the tetrahedron.

\begin{figure}[ht!]
\begin{center}
\begin{tabular}{c c}
\subfigure[$\cT$]{\includegraphics[scale=0.6]{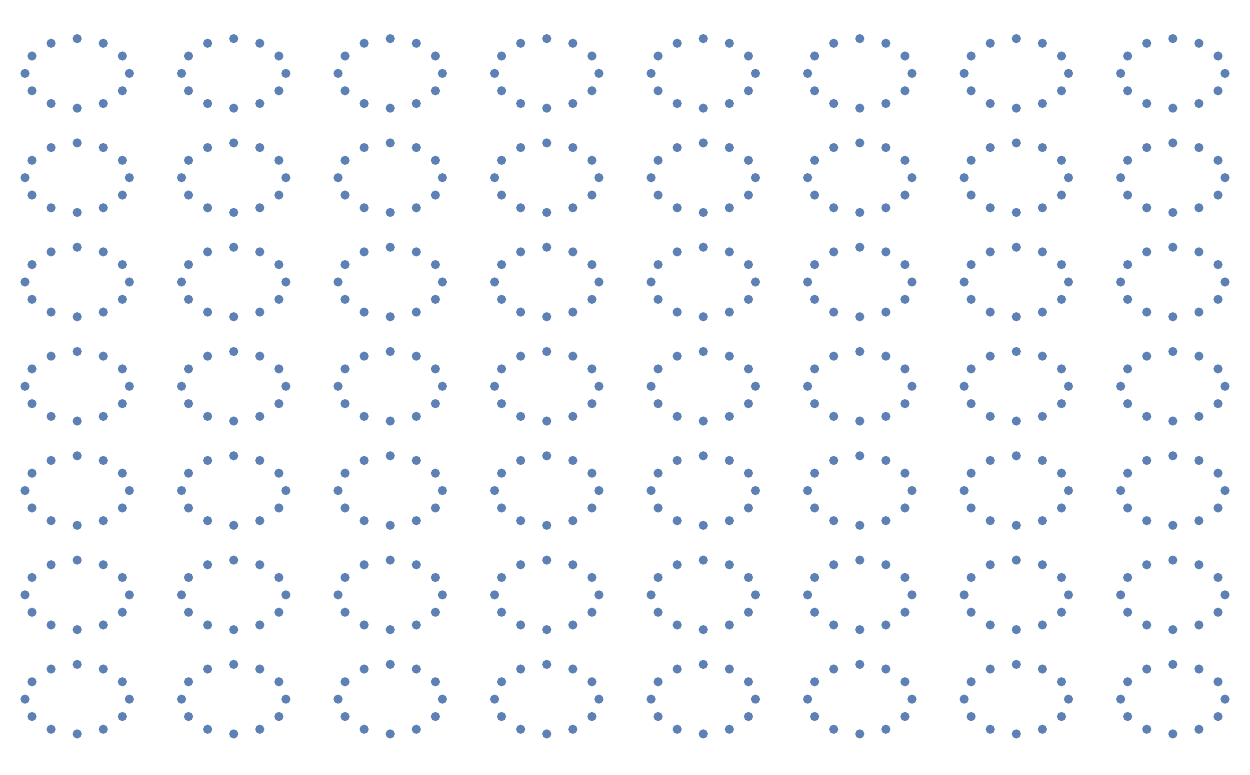}}
&
\subfigure[$A$ acting on $\cT$]{\includegraphics[scale=0.6]{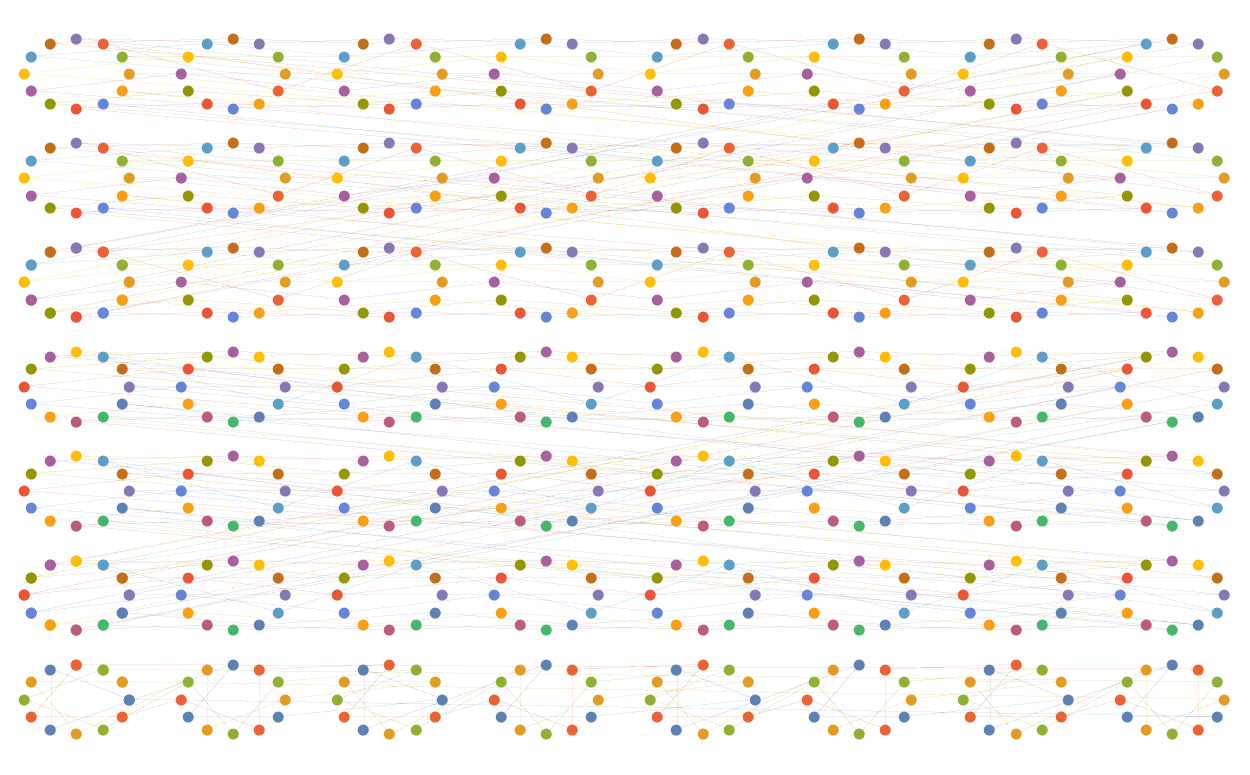}}
\\
\subfigure[$R$ acting on $\cT$]{\includegraphics[scale=0.6]{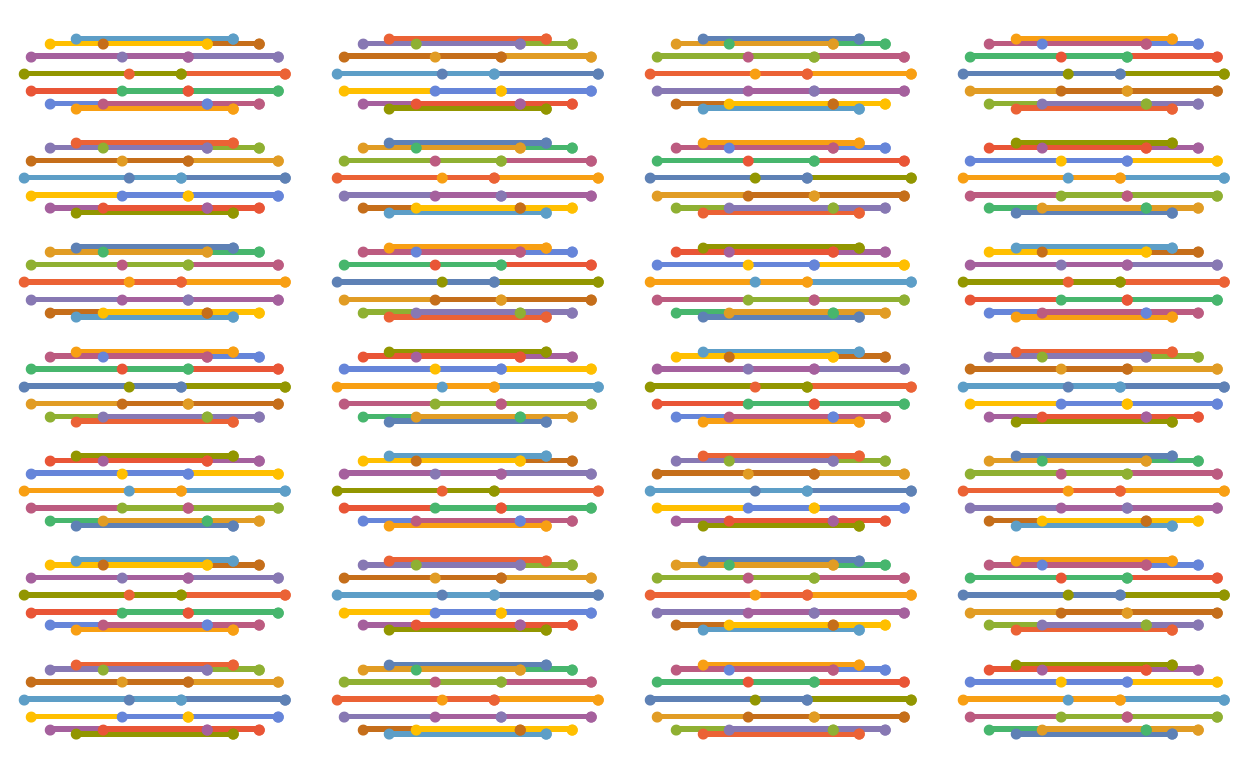}}
&
\subfigure[$S$ acting on $\cT$]{\includegraphics[scale=0.6]{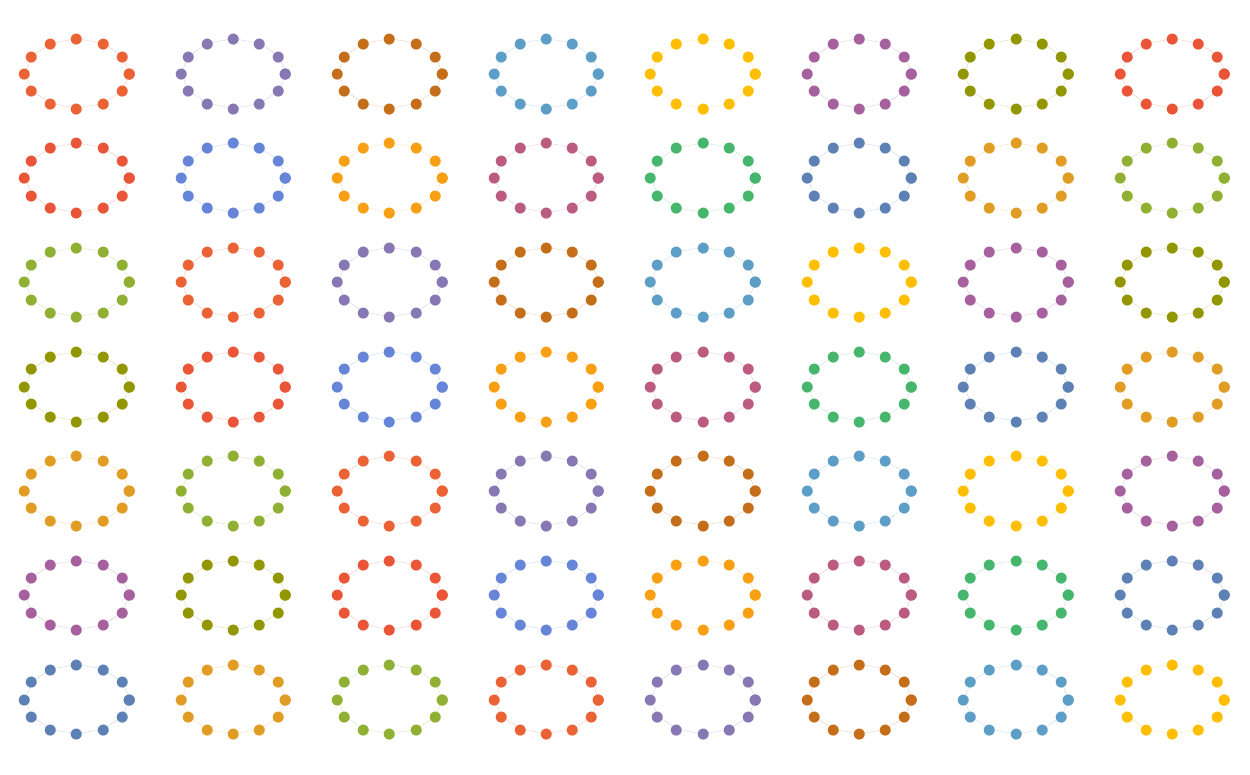}}
\\
\end{tabular}
\end{center}
\caption{Graphical presentation of $A$, $R$, and $S$ acting on the set $\cT$ of all $672$ strong traces of a tetrahedron. Strong traces are presented as vertices of a graph (which consequently has $672$ vertices), two being adjacent when at least one element of $A$ or $R$ or $S$ map one into another. To make the presentation a bit more transparent some edges are left out at figures $(b)$ ad $(d)$. Figure $(b)$ shows $28$ instances of $C_{24}$ which should be replaced with $28$ instances of $K_{24}$, while figure $(d)$ shows $56$ instances of $C_{12}$ which should be replaced with $56$ instances of $K_{12}$.}
\label{fig:acting}
\end{figure}

In Figure~\ref{fig:acting} we graphically present the action of $\Gamma$ on $\cT$ in the case of the tetrahedron. The vertices of a graph on each subfigure represent all $672$ different strong traces of tetrahedron (generated with simple backtracking without eliminating the non-canonical traces). Two vertices $t_1$ and $t_2$ are then adjacent if they lie in the same orbit of $\Gamma$. Note that $\Gamma$ partitions $\cT$ into $3$ orbits of orders $288$, $288$, and $96$. This fact coincides with the results of Table~\ref{tab:platonic}. Subgroups $A$, $R$, and $S$ partition $\cT$ into $28$ orbits of order $24$, $336$ orbits of order $2$, and $56$ orbits of order $12$, respectively.

This is (to our knowledge) the first analyze of the automorphism group of a double trace. We proceed as follows. In 
Section~\ref{sec:branch}, we use the automorphism group
to devise a branch-and-bound 
algorithm that outputs each canonical strong 
double trace of $G$ precisely once. The main idea of the algorithm is avoiding isomorphs by extending minimal forms. Such an idea was first presented in~\cite{rea-1978} where it was called the orderly generation. It is not difficult to see that with minor adjustments, this algorithm can enumerate other varieties of double traces, such as $d$-stable traces, parallel double traces, or antiparallel double traces.
We conclude, in Section~\ref{sec:conclusion}, with some numerical results that 
reveal possible varieties in designing polyhedral polypeptide nanostructures. 

\section{Enumerating strong traces with branch-and-bound strategy}
\label{sec:branch}

In this section we assume that the $n$ vertices of some arbitrary, but fixed, connected graph $G$ with $m$ edges are linearly ordered as $v_0<v_1<\ldots<v_{n-1}$, and that $v_0$, $v_1$ are adjacent. Therefore every cannonical double trace of $G$ starts with $v_0 v_1$. We denote the automorphism group of double traces in graph $G$ with $\Gamma$. To make the arguments more transparent, let $W$ and $W'$ from now on be two different double traces. We first give some additional observations.

\begin{definition}
\label{def:init}
Let $W = w_0 \ldots w_{2m}$ be a double trace of a graph $G$. An {\em initial segment $\init(W)$} of $W$ is the shortest  continuous subsequence of $W$ such that $\init(W)$ starts in $w_0$ and contains all the vertices from $V(G)$.
\end{definition}

\begin{definition}
Let $W = w_0 \ldots w_{2m}$ be a double trace. Then an {\em $i$-initial segment} of $W$, denoted $W_i$, is a subsequence of first $i$ vertices in $W$, i.e., $W_i = w_0 \ldots w_{i-1}$.
\end{definition}

\begin{definition}
A double trace $W$ is {\em $i$-canonical} if and only if for every $\pi \in \Gamma$, the relation $W_i \leq_{\mathit{lex}} \pi(W_i)$ holds.
\end{definition}

\begin{lemma}
If a double trace $W$ of length $2m$ is canonical, then it follows that $W$ is $i$-canonical for all $i$, $1 \leq i \leq 2m$.
\end{lemma}

\proof
Let $W$ be a canonical double trace of length $2m$. Suppose that for some $i$, $1 \leq i < 2m$, $W$ is not $i$-canonical. Then there exists $\pi \in \Gamma$, such that $\pi(W_i) <_{\mathit{lex}} W_i$. Because $W$ is canonical, it follows that $W \leq_{\mathit{lex}} \gamma(W)$ for every $\gamma \in \Gamma$. Therefore, $W \leq_{\mathit{lex}} \pi(W)$. By Definition~\ref{def:lex}, it follows that at the first index $j$, where $w_j \neq \pi(w_j)$, $w_j < \pi(w_j)$. For every $i < j$, $W_i = \pi(W_i)$, while for every every $i \geq j$, $W_i$ contains $w_j$ and therefore $W_i <_{\mathit{lex}} \pi(W_i)$, a contradiction.
\qed

We first explain the auxiliary algorithms used in the main Algoritm~\ref{alg:branch}. If $G$ is a graph with $m$ edges and $p \leq 2m$, then vertex sequence $W_p = w_0 \ldots w_{p-1}$ is a {\em partial double trace} if there exists a double trace $W$ of $G$ for which $W_p$ is its $p$-initial segment. Analogously we define partial double trace for other varieties of double traces (strong and $d$-stable traces). Let $W_p$ be a partial double trace of length $p$. Set $\mathcal{W}$ represent all double traces in $G$ for which their $p$-initial segment is equal to $W_p$. While we say that $W_p$ is lexicographically smaller than different partial double trace $W_p'$ if $W_p=W_p'$ or the first $w_i$, which is different from $w_i'$, is smaller than $w_i'$, we say that $W_p$ is canonical if at least one the double trace from $\mathcal{W}$ is canonical. Stabilizer of a partial double trace $W_p$ is defined as subset of all automorphisms in $\Gamma$ which map at least one double trace from $\mathcal{W}$ back to (not necessary the same) double trace from $\mathcal{W}$. {\em Feasible neighbors} of $w_{p-1}$ in a partial double trace $W_P$ is a subset of its neighbors $N(W_{p-1})$. For every feasible neighbor $v$ then $W_{p+1} = w_0 \ldots w_{p-1} v$ obtained from $W_p$ by adding $v$ also $W_{p+1}$ should be a partial double trace. Analogously for partial strong traces and $d$-stable traces where we have to be careful that $v$ does not cause any new nontrivial repetition of excessive order. For antiparallel or parallel double traces we additionally forbid vertices causing parallel or antiparallel edges in partial double trace, respectively.

Algorithm~\ref{alg:extend} loops through all the feasible neighbors of the last vertex $w_{p-1}$ in a partial double trace $W_p = w_0 \ldots w_{p-1}$ and check which of them, if added to $W_p$ (and obtaining partial double trace $W_{p+1}$), will maintain a canonical partial double trace. Partial double traces obtained in this procedure are added to queue $Q$.

\begin{algorithm*}[h!]
\caption{\textsc{Extend Feasibly}}
\label{alg:extend}
\begin{algorithmic}
\item[] \hspace*{-5mm} {\bf Input}: a partial double trace $W_p = w_0 \ldots w_{p-1}$, $A \subseteq \Gamma$, a queue of $Q$ partial double traces
\STATE{$V' =$ \textsc{Feasible Neighbors}($w_{p-1}$)}
\STATE{$V'' =$ \textsc{Canonical Extension}($V',W_p,A$)}
\FOR {$v \in V''$}
  \STATE {$W_{p+1} = w_1 \ldots w_{p-1}, v$}
  \STATE {$A_v = $ \textsc{Prune}($A,W_{p+1}$)}
  \IF {$W_{p+1}$ is canonical partial double trace}
    \STATE {append $(W_{p+1}, A_v)$ to $Q$}
  \ENDIF
\ENDFOR
\end{algorithmic}
\end{algorithm*}

At each step we use the automorphism group of double traces $\Gamma$ in order to eliminate all partial double traces that would not lead to a construction of a canonical double trace. We achieve that by considering only the lexicographically smallest representative of each orbit of the automorphism group. Simultaneously, we would like to fix vertices that are already in a partial double trace, since we have already checked it for canonicity. Therefore, Algorithm~\ref{alg:prune} returns the automorphisms that are in the stabilizer of partial double trace $W_{p+1}$ (in each step only the last position $p$ has to be checked). Note that until double trace is not completed, we can not determine if automorphism lies in a stabilizer of partial double trace, for all automorphisms from $\Gamma$. Problem is in shifting, since we can not always deteremine how all first $p$ places of shifted partial double trace look like. Therefore we do not discard such automorphisms at this point.

\begin{algorithm*}[h!]
\caption{\textsc{Prune}}
\label{alg:prune}
\begin{algorithmic}
\item[] \hspace*{-5mm} {\bf Input}: set of automorphisms of double traces $A$, partial double trace $W_p$
\item[] \hspace*{-5mm} {\bf Output}: pruned set of automorphisms of double traces $A'$
\STATE{$A' = \emptyset$}
\FOR {$\pi \in A$}
  \IF {$\pi$ in stabilizer of partial double trace $W_p$ or if it can not be determined if $\pi$ is in a stabilizer of $W_p$}
    \STATE {append $\pi$ to $A'$}
  \ENDIF
\ENDFOR
\RETURN $A'$
\end{algorithmic}
\end{algorithm*}

Algorithm~\ref{alg:canonize} loops through all the feasible neighbors of the last vertex $w_{p-1}$ last added to a partial double trace $W_p = w_0 \ldots w_{p-1}$ and denoted with $V \subseteq N(w_{p-1})$. For every $v \in V$ it constructs new partial double trace $W_{p+1} = w_0 \ldots w_{p-1} v$ and analyze orbits of $Aut(G) \cap A$ (no shifts and reverses are allowed, therefore each orbit contains even smaller number of partial double trace) acting on set of these new partial double traces. Then for every such orbit $O$ algorithm select vertex $v \in V$ for which partial double trace $W_p = w_0, \ldots, w_{p-1}$ is lexicographically smallest of partial double traces in $O$. Note that in practice algorithm should only check the position $p$ since Algorithm~\ref{alg:prune} ensures that for every $\pi \in Aut(G) \cap A$ vertices $w_0 \dots w_{p-1}$ are fixed.

\begin{algorithm*}[h!]
\caption{\textsc{Canonical Extension}}
\label{alg:canonize}
\begin{algorithmic}
\item[] \hspace*{-5mm} {\bf Input}: partial double trace $W_p = w_0, \ldots, w_{p-1}$, set of feasible neighbors $V \subseteq N(w_{p-1})$, set of automorphisms $A \subseteq \Gamma$
\item[] \hspace*{-5mm} {\bf Output}: set $V' \subseteq V$ containing for each orbit $O$ of $Aut(G) \cap A$ vertex $v$ for which $W_{p+1} = w_0, \ldots, w_{p-1}, v$ is lexicographically smallest partial double trace of $O$
\IF {$A = \emptyset$ or $A = \{ \text{id} \}$}
  \STATE {return $V$}
\ENDIF
\STATE{$V' = \emptyset$}
\FOR {$v \in V$}
  \STATE {$V'' = \{ v \}$}
  \STATE {$v' = v$}
  \FOR {$\pi \in Aut(G) \cap A$}
    \STATE {append $\pi(v)$ v $V''$}
    \IF {$(w_0, \ldots, w_{p-1}, \pi(v)) <_{\text{lex}} (w_0, \ldots, w_{p-1}, v')$}
      \STATE {$v' = \pi(v)$}
    \ENDIF
  \ENDFOR
  \STATE {append $(v')$ to  $V'$}
  \FOR {$v'' \in V''$}
    \STATE {remove $v''$ from $V$}
  \ENDFOR
\ENDFOR
\RETURN {$V'$}
\end{algorithmic}
\end{algorithm*}

We now present the main Algorithm~\ref{alg:branch}, which enumerates strong traces for an arbitrary graph.

\begin{algorithm*}[h!]
\caption{\textsc{Enumerate Strong Traces}}
\label{alg:branch}
\begin{algorithmic}
\item[] \hspace*{-5mm} {\bf Input}: a graph $G$ with $m$ edges, automorphism group $\Gamma$ of double traces of $G$
\item[] \hspace*{-5mm} {\bf Output}: a list of all non-equivalent double traces $L$
\STATE{$W_1 = v_0 v_1$}
\STATE{$A = \text{Aut}(G)$}
\STATE{$A =$ \textsc{Prune}($A, W_1$)}
\STATE{$Q = \{(W_1, A)\}$}
\WHILE {$Q$ not empty}
  \STATE {$(W, A) =$ head of $Q$}
  \STATE {remove $(W,A)$ from $Q$}
  \IF {$|W| = 2m$}
    \STATE {add $W$ to $L$}
  \ELSE
    \STATE{\textsc{Extend Feasibly}($W, A, Q$)}
  \ENDIF
\ENDWHILE
\RETURN $L$
\end{algorithmic}
\end{algorithm*}

In the rest of the section, we prove the correctness of Algorithm~\ref{alg:branch}.

\begin{theorem}
\label{thm:all-are-cannonical}
Let $W$ be a double trace, which was given as an output of Algorithm~\ref{alg:branch}. Then $W$ is canonical.
\end{theorem}

\proof
Let $W = w_0, \ldots, w_{2m}$ be a double trace obtained as an output of Algorithm~\ref{alg:branch}. Suppose that $W$ is not canonical. Then there exists a double trace $W'$ and $\pi \in \Gamma$,  such that $W' = \pi(W)$ and $\pi(W) <_{\mathit{lex}} W$. Let $i$ be the smallest integer such that $w'_i \neq w_i$. Then $w'_i < w_i$ and $w_j  = w_j' = \pi(w_j)$, for $0 \leq j < i$. For every $1 \leq j < i$ automorphism $\pi$ fixes edge $w_{j-1}w_j$: $w_{j-1}w_j = \pi(w_{j-1})\pi(w_j)$, hence $\pi$ is contained in the stabilizer of $W_j$. Consequently Algorithm~\ref{alg:branch} ~\ref{alg:branch} (Algorithm~\ref{alg:prune} to be more precise) does not eliminate $\pi$ while pruning. Selecting $w_i$ in Algorithm~\ref{alg:branch} (Algorithm~\ref{alg:extend} to be more precise) was not optimal since $w_i'$ would produce lexicographically smaller equivalent double trace. This contradicts the fact that Algorithm~\ref{alg:extend} for every orbit select lexicographically smallest feasible neighbor.
\qed

\begin{theorem}
\label{thm:cannonical-are-all}
Let $W$ be a canonical double trace. Then $W$ is given as an output of Algorithm~\ref{alg:branch}.
\end{theorem}

\proof
Suppose the contrary. Let $W = w_0, \ldots, w_{2m}$ be a canonical double trace which is not given as an output of Algorithm~\ref{alg:branch}. By observations made in Section~\ref{sec:intro}, $W$ starts with $v_0 v_1$. There exists the largest integer $i$ (at least $1$ if no other) such that $W_i$ is the $i$-initial segment of some canonical double trace which is an output of Algorithm~\ref{alg:branch}. Let $\mathcal{W}$ be the set of all (canonical) double traces which are given as an output of Algorithm~\ref{alg:branch} and have $W_i$ as their $i$-initial segment. Let $V_{\mathcal{W},i+1}$ be the set of vertices that lie at the $(i+1)$-th position in (canonical) double traces from $\mathcal{W}$. It follows that $w_{i+1} \notin V_{\mathcal{W},i+1}$. Since $W$ is a double trace, $w_{i+1}$ was in the Algorithm~\ref{alg:branch} (Algorithm~\ref{alg:canonize} to be more precise) part of feasible neighbors of $w_i$ for every double trace from $\mathcal{W}$. Since it was never added it follows that in the same orbit of $Aut(G) \subseteq \Gamma$ than $W$ also lies another (lexicographically smaller) double trace $W' \in mathcal{W}$. That contradicts the fact that $W$ is canonical.
\qed

We presented an algorithm which enumerates all non-equivalent double traces of graph. To enumerate only stron traces all only $d$-stable traces of graph we just have to complement the definition of feasible neighbors.
Da namesto dvojnih obhodov v splo\v snem, \v stejemo le stroge ali le $d$-stabilne obhode je potrebno le spremeniti, kaj so dopustni sosedi vozli\v s\v ca $v$, ki smo ga na $i$-tem koraku dodali v dvojni obhod. Pri strogih obhodih moramo tako paziti, da se ne pojavi nobena netrivialna ponovitev, pri $d$-stabilnih obhodih pa, da se ne pojavi nobena ponovitev reda $\leq d$. Podobno lahko pre\v stevamo tudi paralelne ali antiparalelne dvojne obhode. Za povezavo $e = uv$, ki je bila v dvojnem obhodu \v ze pre\v ckana, si je potrebno zapomniti ali smo jo pre\v ckali v smeri od $u$ proti $v$ ali pa v obratni smeri. Glede na to (in desjtvo ali \v stejemo paralelne ali antiparalelne obhode) potem ustrezno popravimo dopustne sosede vozli\v s\v ca $u$ in $v$.

\section{Concluding remarks and numerical results}
\label{sec:conclusion}

We conclude with some numerical results. In Tables~\ref{tab:platonic},~\ref{tab:prisms}, and~\ref{tab:others} we present enumerations of non-equivalent strong traces for platonic solids, prisms, and some other interesting solids which could be the next candidates to be constructed from coiled-coil-forming segments. Note that $d$, $n$, $m$, $ST$, $aST$, and $pST$ stand for the degree of the graph (if a graph is regular), the number of its vertices and edges, the number and the CPU time used to enumerate strong traces, the number and the CPU time used to enumerate antiparallel strong traces, and number and CPU time used to enumerate parallel strong traces in Tables~\ref{tab:platonic},~\ref{tab:prisms}, and~\ref{tab:others}, respectively. Note that the listed CPU times are measured in seconds. In addition to the number of strong traces, the algorithm for every strong trace also returns its vertex sequence. Therefore it can be used for a thoroughly analysis of some properties that nanostructures self-assembled from these strong traces would have. Further, this analysis help to select a strong trace with the maximal probability to construct a stable nanostructure of desired shape.

\begin{table}[ht!]
\begin{center}
\begin{tabular}{| c | c | c | c | c | c | c | c |}
\hline
graph & d & n & m & \multicolumn{2}{ c |}{ ST } & \multicolumn{2}{ c |}{ pST } \\ \hline 
 & & & & $\#$ & CPU time & $\#$ & CPU time \\ \hline
tetrahedron & 3 & 4 & 6 & 3 & 0.005 & 0 & - \\ \hline
cube & 3 & 8 & 12 & 40 & 0.01 & 0 & - \\ \hline
octahedron & 4 & 6 & 12 & 21479 & 1.86 & 262 & 0.056 \\ \hline
dodecahedron & 3 & 20 & 30 & 2532008 & 2242.31 & 0 & - \\ \hline
\end{tabular}
\end{center}
\caption{Number of strong traces and parallel strong traces for platonic solids}
\label{tab:platonic}
\end{table}

\begin{table}[ht!]
\begin{center}
\begin{tabular}{| c | c | c | c | c | c | c | c |}
\hline
graph & d & n & m & \multicolumn{2}{ c |}{ ST } & \multicolumn{2}{ c |}{ aST } \\ \hline 
 & & & & $\#$ & CPU time & $\#$ & CPU time \\ \hline
$Y_3$ & 3 & 6 & 9 & 25 & 0.007 & 2 & 0.005 \\ \hline
$Y_4$ & 3 & 8 & 12 & 40 & 0.01 & 0 & - \\ \hline
$Y_5$ & 3 & 10 & 15 & 634 & 0.066 & 10 & 0.006 \\ \hline
$Y_6$ & 3 & 12 & 18 & 3604 & 0.377 & 0 & - \\ \hline
$Y_7$ & 3 & 14 & 21 & 21925 & 3.51 & 76 & 0.024 \\ \hline
$Y_8$ & 3 & 16 & 24 & 134008 & 32.5 & 0 & - \\ \hline
$Y_9$ & 3 & 18 & 27 & 833685 & 233.7 & 536 & 0.430 \\ \hline
$Y_{10}$ & 3 & 20 & 30 & 5212520 & 2280.06 & 0 & - \\
\hline
\end{tabular}
\end{center}
\caption{Number of strong traces and antiparallel strong traces for prisms}
\label{tab:prisms}
\end{table}

\begin{table}[ht!]
\begin{center}
\begin{tabular}{| c | c | c | c | c | c | c | c |}
\hline
graph & d & n & m & \multicolumn{2}{ c |}{ ST } & \multicolumn{2}{ c |}{ aST } \\ \hline 
 & & & & $\#$ & CPU time & $\#$ & CPU time \\ \hline
$4$-pyramid & - & 5 & 8 & 52 & 0.004 & 4 & 0.008 \\ \hline
$3$-bipyramid & - & 5 & 9 & 470 & 0.013 & 0 & - \\ 
\hline
\end{tabular}
\end{center}
\caption{Number of strong traces and antiparallel strong traces in $4$-pyramid and $3$-bipyramid}
\label{tab:others}
\end{table}

All the calculations were made using Algorithm~\ref{alg:branch} and computational resources at SageMathCloud~\cite{sage}. It was observed in~\cite{fi-2013}, that a graph $G$ admits a parallel strong trace if and only if $G$ is Eulerian, and that $G$ admits an antiparallel strong trace if and only if there exists a spanning tree $T$ of $G$ with the property that every component of the co-tree $G-E(T)$ is even. Therefore, we omit the information about antiparallel and parallel strong traces for graphs not admitting them. Some of these calculations were already presented in~\cite{gr-2013} and~\cite{kl-2013}.

Another possible approach to strong trace construction exploits the observation that a strong trace can be nicely drawn on a surface in which the given graph is embedded. This surface can be cut along certain edges which results in one or more surfaces with boundary. Each of the resulting surfaces with boundary carries a part of the information about the strong trace. The strong trace can be reconstructed by gluing those smaller pieces back together. This topological approach will be elaborated in~\cite{bas-2016b}.

\section*{Acknowledgments}

The authors would like to thank to Anders Skovgaard Knudsen who independently calculated the number of strong traces in platonic solids and shared the results for comparison.

This research was supported in part by Slovenian Research Agency under research grants L7-5554 and P1-0294.



\begin{thebibliography}{10}

\bibitem{bas-2016b}
  N.~Ba\v si\' c, D.~Bokal, T.~Pisanski, J.~Rus,
  \emph{Graph embeddings yield natural strong trace realizations}, in preparation.

\bibitem{ben-1998}
E.~Benevant~L\' opez, D.~Soler~Fern\' andez, \emph{Searching for a strong double tracing in a graph},
  Sociedad de Estad\' istica e Investigaci\' on Operativa Top Vol. 6 (1998), 123--138.

\bibitem{br-1988}
H.~J. Broersma and F.~G{\"o}bel, \emph{{$k$}-{T}raversable graphs}, Ars Combin.
  \textbf{29} (1990), no.~A, 141--153, Twelfth British Combinatorial Conference
  (Norwich, 1989).

\bibitem{sage}
The~Sage Developers, \emph{{S}age {M}athematics {S}oftware ({V}ersion
  6.9.(2015-10-10)}, 2015, {\tt http://www.sagemath.org}.

\bibitem{eg-1984}
R.~B. Eggleton and D.~K. Skilton, \emph{Double tracings of graphs}, Ars Combin.
  \textbf{17} (1984), no.~A, 307--323.

\bibitem{el-2013}
J.~A. Ellis-Monaghan, A.~McDowell, I.~Moffatt, and G.~Pangborn, \emph{D{NA}
  origami and the complexity of {E}ulerian circuits with turning costs}, Nat.
  Comput. \textbf{14} (2015), no.~3, 491--503.

\bibitem{eul-1736}
L.~Euler, \emph{Solutio problematis ad geometriam situs pertinentis},
  Commentarii Academiae Scientiarum Imperialis Petropolitanae \textbf{8}
  (1741), 128--140.

\bibitem{fi-2013}
G.~Fijav{\v{z}}, T.~Pisanski, and J.~Rus, \emph{Strong traces model of
  self-assembly polypeptide structures}, MATCH Commun. Math. Comput. Chem.
  \textbf{71} (2014), no.~1, 199--212.

\bibitem{fl-1990}
H.~Fleischner, \emph{Eulerian graphs and related topics. {P}art 1. {V}ol. 1},
  Annals of Discrete Mathematics, vol.~45, North-Holland Publishing Co.,
  Amsterdam, 1990.

\bibitem{fl-1991}
H.~Fleischner, \emph{Eulerian graphs and related topics. {P}art 1. {V}ol. 2},
  Annals of Discrete Mathematics, vol.~50, North-Holland Publishing Co.,
  Amsterdam, 1991.

\bibitem{fur-1988}
M.~L. Furst, J.~L. Gross, and L.~A. McGeoch, \emph{Finding a maximum-genus
  graph imbedding}, J. Assoc. Comput. Mach. \textbf{35} (1988), no.~3,
  523--534.
  
\bibitem{gab-1986}
  H.~N.~Gabow, M.~Stallman, \emph{An Augmenting Path Algorithm for Linear Matroid Parity},
  Combinatorica 6 2 (1986), 123-150.

\bibitem{god-2001}
C.~Godsil and G.~Royle, \emph{{A}lgebraic {G}raph {T}heory}, Graduate Texts in
  Mathematics, vol. 207, Springer-Verlag, New York, 2001.

\bibitem{gr-2013}
H.~Gradisar, S.~Bo{\v{z}}i{\v{c}}, T.~Doles, D.~Vengust,
  I.~Hafner~Bratkovi{\v{c}}, A.~Mertelj, B.~Webb, A.~{\v{S}}ali,
  S.~Klav{\v{z}}ar, and R.~Jerala, \emph{Design of a single-chain polypeptide
  tetrahedron assembled from coiled-coil segments}, Nat. Chem. Biol. \textbf{9}
  (2013), no.~6, 362--366.

\bibitem{kl-2013}
S.~Klav{\v{z}}ar and J.~Rus, \emph{Stable traces as a model for self-assembly
  of polypeptide nanoscale polyhedrons}, MATCH Commun. Math. Comput. Chem.
  \textbf{70} (2013), no.~1, 317--330.

\bibitem{ko-2015}
V.~Ko{\v{c}}ar, S.~Bo{\v{z}}i{\v{c}}~Abram, T.~Doles, N.~Ba{\v{s}}i{\'{c}},
  H.~Gradi{\v{s}}ar, T.~Pisanski, and R.~Jerala, \emph{Topofold, the designed
  modular biomolecular folds: polypeptide-based molecular origami
  nanostructures following the footsteps of dna}, WIREs Nanomed.
  Nanobiotechnol. \textbf{7} (2015), no.~2, 218--237.

\bibitem{pis-2009}
T.~Pisanski and A.~{\v{Z}}itnik, \emph{Representing graphs and maps}, Topics in
  topological graph theory, Encyclopedia Math. Appl., vol. 128, Cambridge Univ.
  Press, Cambridge, 2009, pp.~151--180.

\bibitem{rea-1978}
R.~C.~Read, \emph{Every one a winner},
  Annals Discrete Math. 2 (1978), 107--120.

\bibitem{rus-2015}
J.~Rus, \emph{Antiparallel $d$-stable traces and a stronger version of {O}re
  problem},  (2015), submitted.

\bibitem{sa-1977}
G.~Sabidussi, \emph{Tracing graphs without backtracking}, Methods of Operations
  Research {XXV}, Part 1 (R.~Henn, P.~Kall, B.~Korte, O.~Krafft, W.~Oettli,
  K.~Ritter, J.~Rosenm{\"{u}}ller, N.~Schmitz, H.~Schubert, and W.~Vogel,
  eds.), First Symposium on Operations Research, University of Heidelberg,
  Verlag Anton Hain, June 1977, pp.~314--332.

\bibitem{tarry-1895}
G.~Tarry, \emph{Le probl\`eme des labyrinthes}, Nouv. Ann. \textbf{3} (1895),
  no.~XIV, 187--190.

\bibitem{th-1990}
C.~Thomassen, \emph{Bidirectional retracting-free double tracings and upper
  embeddability of graphs}, J. Combin. Theory Ser. B \textbf{50} (1990),
  198--207.

\bibitem{we-1996}
D.~B. West, \emph{{I}ntroduction to {G}raph {T}heory}, Prentice Hall, Upper
  Saddle River, 1996.

\end{thebibliography}
\end{document}